\documentclass{conm-p-l}
 \usepackage{amsthm,amsmath,amsfonts,amssymb,amscd}

 \usepackage[all]{xy}

 \CompileMatrices

\newtheorem{Theorem}{Theorem}[section]
\theoremstyle{definition}
\newtheorem{Definition}[Theorem]{Definition}

\newtheorem{Lemma}[Theorem]{Lemma}

\theoremstyle{remark}
\newtheorem{Remark}[Theorem]{Remark}

\newtheorem{Example}[Theorem]{Example}

\numberwithin{equation}{section}

\begin{document}

\title
 {A Counterexample to a conjecture of Bosio and Meersseman}

\author{ David Allen }
\address{Department of Mathematics Iona College, New Rochelle, 10801}
\email{dallen@iona.edu}
\author{Jos\'{e} La Luz}
\address{Department of Mathematics, University of Puerto Rico,
PO BOX 23355 San Juan, Puerto Rico 00931-3355}
\email{jluz@cnnet.upr.edu } 


\maketitle

\begin{abstract}

In \cite{BM} the following is conjectured:  If $P$ is dual
neighborly, then $Z_P$ is diffeomorphic to the connected sum of
products of spheres. In this paper a counterexample is provided.  \\

\end{abstract}

\section*{Introduction}

     In \cite{A} methods from unstable homotopy theory were
     used to extend homotopy computations for various spaces that
     are studied in Toric Topology. In that paper, one of the main
     theorems shows that through a range the higher homotopy groups
     of the Borel space are isomorphic to the higher
     homotopy groups of a certain wedge of spheres.  The range
     restriction arises from the combinatorics of $P$ and depends
     on the monomials that appear in the quotient in the face ring
     $\mathbb{Z}(P)$.  The calculations in that paper provide
     enough homotopy theoretic information that allows for one to
     compare $\pi_*(Z_P)$ with $\pi_*(B_TP)$ via the fibration
     $Z_P \rightarrow B_TP \rightarrow BT^m$ to show that the
     conjecture can not hold.  The computations are rather
     straightfoward, however, they illustrate a particularly
     useful application of \ref{borel homotopy} when it comes to
     answering or at least testing problems of this type.
     In Toric Topology it is well known by
     the work of Buchstaber, Panov and their collaborators that
     the higher homotopy groups of the moment angle complex $Z_P$
     are isomorphic to the higher homotopy groups of the
     complement of a complex coordinate subspace arrangement,
     sometimes denoted in the literature as $U(K)$.  In \cite{BM} the
     "link" $X$ that appears in their conjecture has the same higher homotopy
     groups as the spaces that appear in Toric topology, in
     particular, the moment angle complex $Z_P$ as well as the
     Borel space $B_TP$.  This is pointed out by \cite{BM} in the
     proof of 7.6 and 7.7. \\
     \indent The paper is set up as follows.  In \S 1 the
     combinatorics needed in the sequel will be discussed. In \S 2
     the main definitions from Toric Topology will be listed.  In
     \S 3 theorems from homotopy theory will be listed and in \S 4 the main result will be
     proven.

\section{Combinatorics}
   All of the definitions in this section can be found in
   \cite{BP,B}.  For a more in-depth
   discussion of the material in this section the reader is urged
   to refer to one of these references.

\begin{Definition}
    The \textsl{affine hull} of the points $x_1,..,x_n$ where
    $x_i\in\mathbb{R}^n$ is the set $\{\sum_{i=1}^n a_ix_i | a_i\in\mathbb{R},
    a_1 +\cdot\cdot\cdot+ a_n = 1\}.$\\

\end{Definition}
  A \textit{convex polytope} is the convex hull of finitely many points
    in $\mathbb{R}^n$.  We assume that all polytopes $P$ contains 0 in its interior.
    $(\mathbb{R}^n)^*$ will denote the vector space dual to $\mathbb{R}^n$.
    There is an equivalent definition.
\begin{Definition}
    A convex polyhedron $P$ is an intersection of finitely many
    half-spaces in $\mathbb{R}^n$.

$$
\begin{CD}
P = \{\overrightarrow{x}\in \mathbb{R}^n | \langle
            \overrightarrow{l_i}, \overrightarrow{x}\rangle \geq
            -a_i, i = 1,..,m\}\\
\end{CD}
$$
\end{Definition}

    The dimension of $P$ is the dimension of its affine hull.
   A \textit{supporting hyperplane} $H$ of $P$ is an affine hyperplane
    which intersects $P$ such that $P$ is contained in one the
    closed half-spaces determined by $H$.  A \textit{face}
     of $P$ is $P \cap H$ where $H$ is a supporting
    hyperplane. The \textit{boundary} of $P$, denoted $\partial P$ is the union of all
    the proper faces of $P$.  The vertices of $P$ are the 0-dimensional faces.  The
    \textit{facets} are the $(n-1)$-dimensional faces.  That is, the co-dimension
    1-faces.  An $n$-dimensional convex polytope is \textit{simple} if the number of
    facets meeting at each vertex is exactly $n$.  For each co-dimension $k$ face
    $F_k = \bigcap_{j=1}^k F_{\iota_j}$. $P$ is said to be {\it $q$-neighborly} if
   $\bigcap_{j=1}^q F_{\iota_j}\neq \emptyset$ for any $q$ facets.  For example,
   the square is 1 neighborly.   An $n$-dimensional convex polytope is called \textit{simplicial}
   if there are at least n facets meeting at each vertex.

\begin{Definition}
  For any convex polytope $P \subset \mathbb{R}^n$ define its
  polar set $P^* \subset (\mathbb{R}^n)^*$ by
$$
\begin{CD}
P^* = \{ \textbf{x}' \in (\mathbb{R}^n)^* |\; \langle \textbf{x}',
\textbf{x} \rangle \geq -1 \; \forall\; \textbf{x} \in P\}\\
\end{CD}
$$

\end{Definition}

\indent  $P^*$ is a convex polytope since $0 \in P$.  We recall some
properties of cyclic polytopes \cite{BP,B}. For $d \geq 2$ the
moment curve $M_d$ in $\mathbb{R}^d$ is the curve parametrized by
$t\mapsto x(t) = (t,t^2,..,t^d), t \in \mathbb{R}^d$.  A
\textit{cyclic polytope} of type $C(n,d)$ where $n\geq d+1$ and $d
\geq 2$ is the convex hull of $\{x(t_1),..,x(t_n)\}$ where the $t_i$
  are distinct real numbers.

\indent It is well known that a simplicial polytope is {\it k
neighborly} if any $k$ vertices span a face \cite{BP}.  We say a
simplicial $n$ polytope with an arbitrary number of vertices is {\it
neighborly} if it is $[\frac{n}{2}]$ neighborly.  If $P$ is a
neighborly polytope, then its dual $P^*$ is said to be {\it dual
neighborly}. In what follows we need to describe what sets of
vertices span a face of the polytope $C(n,d)$.  We assume that $t_1
< t_2 < ..< t_n$.  Let $\{x(t_1),..,x(t_n)\}$ be the vertex set of
$C(n,d)$.

\begin{Definition}
Let $X \subset \{x(t_1),..,x(t_n)\}$.  A component of $X$ is a
non-empty subset $Y = \{x(t_j), x(t_{j+1}),..,x(t_{k-1}),
x(t_k)\}$ of $X$ such that $x(t_{j-1}) \notin X$ for $j >1$ and
$x(t_{k+1}) \notin X$ for $k <n$.

\end{Definition}

\indent A component $Y$ is {\it proper} if
$x(t_1)$ and $x(t_n)$ do
not belong to the set \newline $\{x(t_j), x(t_{j+1}),..,x(t_{k-1}),
x(t_k)\}$.  A component containing an odd number of vertices is
called an {\it odd component}.  The following theorem from \cite{B}
describes when a set of vertices span a face of a cyclic polytope.

\begin{Theorem}\label{faces of cylcic poly}
Let $C(n,d)$ be a cyclic polytope such that $t_1 < t_2 < ..< t_n$.
Let $X$ be a subset of $\{x(t_j), x(t_{j+1}),..,x(t_{k-1}),
x(t_k)\}$ containing $k$ points where $k \leq d$.  Then $X$ is the
set of vertices of a $(k-1)$ face of $C(n,d)$ if and only if the
number of odd components of $X$ is at most $d-k$.

\end{Theorem}
The relation between cyclic polytopes and neighborliness is given by
the following theorem that appears in \cite{B}.

\begin{Theorem}
  A cyclic polytope of type $C(n,d)$ is a simplicial $k$ polytope
  for $k \leq [\frac{d}{2}]$.
\end{Theorem}

It is stated in \cite{B} pg. 92 that the dual of a cyclic polytope
is a dual neighborly polytope.

In what follows $P$ is an $n$ dimensional, {\it q-neighborly}, {\it
simple convex polytope} unless otherwise stated.

\begin{Definition}
Let $F= \{F_1,..,F_m\}$ be the set of facets of $P$.  For a fixed
commutative ring $R$ with unit we have

$$
\begin{CD}
R(P)
    =R[v_1,..,v_m]/ I = \langle v_{i_1}
\cdot\cdot\cdot v_{i_k}| \bigcap_{j = 1}^k F_{i_j} =
\emptyset\rangle\\
\end{CD}
$$
\newline
where $|v_i| = 2$ are indexed by the facets and the ideal $I$ is
generated by square free monomials.

\end{Definition}
  In the sequel we will be interested in $\mathbb{Z}(P)$.
  We denote by $|I|$ the cardinality of of the generating set of the ideal.
  There is a relation between $\mathbb{Z}(P)$ and
  $\mathbb{Z}(P^*)$.  Given $P$ and its polar $P^*$.
  Let $K_P$ be the boundary of $P^*$, then $\mathbb{Z}(P) \cong \mathbb{Z}(K_P)$.

\section{Toric Spaces and the Borel Space}
       Let $P$ be as in the previous section; $T^m$ is the $m$ dimensional topological
      torus which is a product of circles indexed by the facets of
      $P$ and $BT^m$ the classifying space for $T^m$.   For each co-dimension $k$ face $F_k = \bigcap_{j=1}^k
       F_{\iota_j}$ there is a coordinate subgroup  $T^{F_k} = S_{F_{\iota_1}}^1 \times
S_{F_{\iota_2}}^1\times \cdot\cdot\cdot \times S_{F_{\iota_k}}^1$.

\begin{Definition}\label{Borel Space}
     Let a space $X$ be endowed with a torus action.  Then the
     Borel space $B_TX$ is the identification space

$$
\begin{CD}
ET^m \times X/\sim$ =  $ET^m \times_{T^m} X\\
\end{CD}
$$
    where the equivalence relation is defined by: $(e,x)\sim
    (eg,g^{-1}x)$ for any $e\in ET^m$ and $x\in X$ , $g\in T^m$\\

\end{Definition}

\begin{Remark}

    $B_TP$ refers to $ET^m
    \times X/\sim$ which has the homotopy type as $P$ when $P$ is a
    simple convex polytope \cite{DJ}.  There is no ambiguity in writing
    $B_TP$ instead of $B_TX$-the Borel construction applied to a space.
    It is also shown in \cite{DJ} that $H^*(B_TP) \cong \mathbb{Z}(P)$

\end{Remark}

    Sometimes one refers to this as applying the Borel construction
    to $X$. One nice property of this definition is that it allows for one
    to move the torus action into the quotient.
    Another very important characteristic of the Borel space is the existence of the
    fibration $X \rightarrow  ET^m \times_{T^m} X \rightarrow  BT^m$.

We take as our definition of the \textit{moment angle complex}
\cite{BP}, \cite{DJ}

\begin{Definition}

$$
\begin{CD}
Z_P =  T^m \times P / \sim\\
\end{CD}
$$
\newline

where $(g,p) \sim (h,q) \Leftrightarrow p = q \in F_k$  and
$g^{-1}h\in  T^{F_k}$\\

\end{Definition}

The reader should note that this definition carries over immediately
when $P$ is a simple polyhedral complex.  However, it is not the
case that $Z_P$ is a manifold. $T^m$ acts on $Z_P$ \cite{BP} from
which the existence of the fibration $Z_P \rightarrow B_TP
\rightarrow BT^m$ follows.  In fact, $Z_K$ exists for any simplicial
complex $K$. The Borel construction applied to the corresponding
moment angle complex gives the space $B_TZ_K = ET^m \times_{T^m}
Z_K$ and the associated fibration $Z_K \rightarrow B_TZ_K
\rightarrow BT^m$.

If $K$ is a simpicial sphere then $Z_K$ is a manifold. It should be
noted that the moment angle complex behaves well with respect to
products and other operations on the level of the polytope
\cite{BP}. $U(K)$ is the complex coordinate subspace arrangement
complement associated to $K$.   We now give its definition.  First,
we assume that $K$ is an $n-1$ dimensional simplicial complex on a
vertex set $\{v_1,..,v_m \}$ denoted by $[m]$.  Let $\sigma =
\{\iota_1,..,\iota_k\} \subset [m]$.   A  \textit{coordinate
subspace} of $\mathbb{C}^m$ is $L_\sigma =\{(z_1,..,z_m)\in
\mathbb{C}^m | z_{\iota_1} = \cdot \cdot \cdot = z_{\iota_k} = 0 \}$
\cite{BP}. The dimension of $L_\sigma = m - |\sigma|$. An
\textit{arrangement} $A = \{L_1,..,L_j\}$ is \textit{coordinate} if
each $L_i$ is coordinate.   The \textit{complex coordinate subspace
arrangement} associated to $K$, $CA(K) = \{L_\sigma|\sigma\notin
K\}$ and the \textit{complex coordinate subspace arrangement
complement} is  $U(K) = \mathbb{C}^m \backslash
\bigcup_{\sigma\notin K} L_\sigma$. In \cite{BP} it is shown that
there is a bijection between the set of simplicial complexes on
$[m]$ and the set of coordinate subspace arrangement complements in
$\mathbb{C}^m$. We list two examples taken from \cite{BP}.
\begin{Example}
    If $K = \partial\Delta^{m-1}$ then  $U(K) =
    \mathbb{C}^m\backslash \{0\}$

\end{Example}
and
\begin{Example}
    If $K_P$ is dual to an $m$-gon then
$$
\begin{CD}
    U(K) = \mathbb{C}^m\backslash \bigcup_{i-j\; \neq 0,1\; mod\;(m)}
    \{z_i = z_j =0\}\\
\end{CD}
$$

\end{Example}

\cite{BP} proved that $Z_K \subset U(K)$ and that there is an
equivariant deformation retraction $U(K) \rightarrow Z_K$.  In
particular, any information concerning the homotopy of $Z_K$ gives
information about the homotopy of $U(K)$. It is clear  from this
theorem and the fibration:  $Z_K \rightarrow B_TP \rightarrow BT^m$
that $\pi_*(U(K)) \cong \pi_*(B_TP)$ for $* >2$. This will be
essential in showing that the conjecture can not be true.

\section{Homotopy Theory}
      In this section we recall the Hilton-Milnor theorem,
      one of the main theorems from \cite{A} as well
      as some terminology that appears in the sequel.\\

Let $A_{k}$ be the free nonassociative algebra on
$\{x_{1},x_{2},\dots ,x_{k}\}$. If $w\in A$ is a monomial we define
the weight of $w$, $W(w)$, to be the number of its factors and
$a_{i}(w)$ be the number of times that the generator $x_{i}$ appears
in $w$. Following \cite{W}, we single out certain elements in
$A_{k}$ called basic products. We single out
$\mathcal{A}_{i}(A_{k})$ the basic products of weight $i$. We define
$\mu :\mathbb{N} \to \mathbb{N}$ such that

\[
\mu(n)=
\begin{cases}
1 & n=1\\
(-1)^{k}& n=p_{1}\cdots p_{k},\; p_{i}\ne p_{j}, \text{$p_{i}$ prime}, \;n>1\\
0 & \text{otherwise}
\end{cases}
\]

\begin{Theorem}
The number of elements in  $\mathcal{A}_{n}(A_{k})$ is given by
the formula
\begin{equation}\label{eq}
\frac{1}{n}\sum _{d|n}\mu(d)k^{\frac{n}{d}}
\end{equation}
\end{Theorem}
Let $X_{1}, X_{2}, \ldots,X_{k}$ be connected spaces with
basepoint.  For $w\in \mathcal{A}_{\ast}(A_{k})$ a monomial we
define
$$w(X_{1}, X_{2}, \ldots,X_{k})=(X_{1})^{a_{1}(w)}\wedge (X_{2})^{a_{2}(w)}\wedge \ldots (X_{k})^{a_{k}(w)}$$
Where $(X)^{k}$ is $X$ smashed with itself $k$ times.

\begin{Theorem}[Hilton-Milnor]
There is a map
\begin{equation}\label{hm}
 h: J(X_{1}\vee X_{2}\vee \ldots \vee X_{k})\to \prod _{w\in \mathcal{A}_{\ast}(A_{k})}Jw(X_{1},\ldots,X_{k})
\end{equation}
where $J$ is the James reduced product. Moreover this map is a
homotopy equivalence.
\end{Theorem}

In the case $k=16$ and $X_{i}=S^{4}$ for every $1 \leq i \leq 16$.
The left side of (\ref{hm}) gives
$$J(\vee_{16}S^{4})\simeq \Omega \Sigma (\vee_{16} S^{4})\simeq \Omega(\vee_{16}S^{5})$$
and the right side of (\ref{hm}) gives

\begin{align}\notag
\prod _{w\in \mathcal{A}_{\ast}(A_{16})}Jw(S^{4},\ldots,S^{4})
\notag &\simeq \prod _{w\in \mathcal{A}_{\ast}(A_{16})}J(S^{4
a_{1}(w)+\cdots+4 a_{k}(w)}) \\ \notag
&\simeq \prod _{w\in \mathcal{A}_{\ast}(A_{16})}\Omega \Sigma S^{4(a_{1}(w)+\cdots +a_{k}(w))} \notag \\
&\simeq \Omega \prod _{w\in
\mathcal{A}_{\ast}(A_{16})}S^{(4W(w)+1)}\notag
\end{align}
The number of spheres of dimension $4W(w)+1$ can be calculated using
formula (\ref{eq}). All of this implies
\begin{equation}\label{f}
\pi_{\ast}(\vee_{16}S^{5})\cong \prod_{w\in \mathcal{A}_{\ast}(A_{16}) }
\pi_{\ast}(S^{(4W(w)+1)})
\end{equation}
\indent We briefly summarize the notion of relations among relations
that appears in \cite{A}. Consider the following 2 dimensional
simple convex polytope.

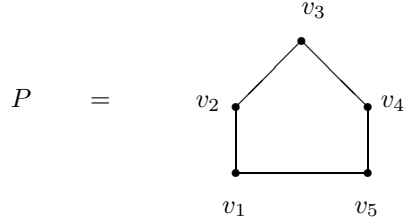
\begin{picture}(50,50)(100,400)
\def\mp{\multiput}
\def\elt{\circle*{3}}

\put(150,375){${P}$}

\put(180,375){${=}$}

 \mp(260,400)(200,320){1}{\line(-1,-1){25}}

\mp(260,400)(200,320){1}{\line(1,-1){25}}

\mp(235,375)(0,30){1}{\line(0,-1){25}}

\mp(285,375)(0,30){1}{\line(0,-1){25}}

\mp(235,350)(0,30){1}{\line(1,0){50}}

\mp(235,375)(0,0){1}{\elt}

 \put(220,375){${v_2}$}

\mp(235,350)(0,0){1}{\elt}

\put(230,335){${v_1}$}

\mp(285,350)(0,0){1}{\elt}

\put(280,335){${v_5}$}

\mp(285,375)(0,0){1}{\elt}

\put(290,375){${v_4}$}

\mp(260,400)(0,0){1}{\elt}

\put(260,410){${v_3}$}\\

\end{picture}

\vspace*{1.1in}

It is clear that the face ring $\mathbb{Z}(P) \cong
\mathbb{Z}[v_1,..,v_5]/ \langle v_1v_3, v_2v_4, v_3v_5, v_4v_1,
v_5v_2\rangle$.  Resolving the ideal in the category of free
$\mathbb{Z}$ algebras produces a resolution:

$$
\begin{CD}
R^* @>\iota^*>> C^* @>p^*>> \mathbb{Z}(P)\\
\end{CD}
$$
\newline

where $p^*$ is the induced map in cohomology coming from $p: B_TP
\rightarrow BT^m$, $R^*$ is the free algebra generated by $ker
\;p^*$ and $C^*$ is the free algebra. More explicitly, for this
example we have

$$
\begin{CD}
\mathbb{Z}[x_1,..,x_5] @>\iota^* >> \mathbb{Z}[v_1,..,v_5] @>>>
\mathbb{Z}(P)\\
\end{CD}
$$
\newline

where $\iota^*$ sends the $x_i$ to the monomials that generate the
ideal $I$.  For example,  $\iota^*(x_1) = v_1v_3$ and
$\iota^*(x_2) = v_2v_4$.

\begin{Definition}\label{Def of relation aong relation}
    Let $\overline{I}$ and $I'$ be two multi-indexes such that $\overline{I} \neq I'$ and
   $\overline{I}, I' \subset S = \{\iota_1,..,\iota_m\}$.  Suppose $i \neq j$.
   We call a relation of the form:

$$
\begin{CD}
\iota^*(x_i)\prod_{k \in \overline{I}} v_k - \iota^*(x_j)\prod_{k' \in I'} v_{k'} = 0\\
\end{CD}
$$

in $C^*$  a relation among relations. We denote such a relation by
$\overline{R}$.

\end{Definition}
In the example above, a relation among relations is  $\overline{R} =
    \iota^*(x_1)v_2v_4
- \iota^*(x_2)v_1v_3$. The degree of this relation among relations
denoted by $| \overline{R}|$ is  8. $\Re_{min}$ is a relation among
relations of the smallest degree. In general, $\Re_{min}$ is not
unique as one can readily verify from the example above. In \cite{A}
the $BP$ analogue of this procedure was carried out with the idea of
setting up the Unstable Adams Novikov Spectral Sequence through a
range.  Let $r_j \in I$.

\begin{Theorem}\label{borel homotopy}\cite{A}
Suppose $P$ is a $q$ neighborly, $n$ dimensional simple polyhedral
complex and $I$ is the ideal in the face ring $\mathbb{Z}(P)$. Let
$p$ be a prime number,  $m$ be the number of facets of $P$. \\

$ \pi_{s-1}(B_TP)_{(p)} = \pi_{s-1}((\bigvee_{r_j \in I } S^{|r_j|-1}))_{(p)} $ for $s \leq |\Re_{min}|- 1$ \\

and $\pi_2(B_TP) = \mathbb{Z}^{\oplus m}$\\

\end{Theorem}

\section{Main Result}

    In \cite{BM}, it is shown that when $P$ is the cyclic ploytope $C(8,4)$ the
    cohomology ring of the moment angle complex is \\
$$
\begin{CD}
H^*(Z_P) \cong H^*((\sharp_{16} S^5 \times S^7) \sharp
(\sharp_{15}
S^6 \times S^6))\\
\end{CD}
$$
\newline
    It was stated in \cite{BM} that it was unknown whether or not
    $Z_P$ is diffeomorphic to the connected sum $(\sharp_{16} S^5 \times S^7) \sharp
(\sharp_{15} S^6 \times S^6)$.  We compute one of the higher
homotopy groups of this connected sum and apply Theorem \ref{borel
homotopy} to compare it to the same higher homotopy group of $B_TP$.
We need a few preliminary results.

\begin{Lemma}\label{facering}
    For $P = C(8,4)$,  the face ring $\mathbb{Z}(P) = \mathbb{Z}[v_1,..,v_8] /
    I$ where $I$ is generated by the square free monomials:

\end{Lemma}

\begin{itemize}
  \item $v_1v_3v_5$,  $v_1v_3v_6$,  $v_1v_3v_7$,  $v_1v_4v_6$
  \item $v_1v_4v_7$,  $v_1v_5v_7$,  $v_2v_4v_6$,  $v_2v_4v_7$
  \item $v_2v_4v_8$,  $v_2v_5v_7$,  $v_2v_5v_8$,  $v_2v_6v_8$
   \item $v_3v_5v_7$, $v_3v_5v_8$,  $v_3v_6v_8$,  $v_4v_6v_8$
\end{itemize}

\begin{proof}
   This follows from an application of Theorem \ref{faces of cylcic poly} applied to the cyclic polytope $P = C(4,8)$.

\end{proof}

\begin{Remark}

From this explicit description of the face ring one readily sees
that $|\Re_{min}| = 8$.  It is necessary to know what the range is
when applying \ref{borel homotopy}.  The reader should note that
the input for \ref{borel homotopy} consists of $|\Re_{min}|$,
which in many cases depends on computing the face ring
$\mathbb{Z}(P)$ explicitly.
\end{Remark}

We compute the integral homology of $ (\sharp_{16} S^5 \times S^7)
\sharp
(\sharp_{15} S^6 \times S^6)$.  Let $T_{m,n}=S^{m}\times S^{n}$ \\

\begin{Lemma}

Let $M= (\sharp_{16} T_{5,7}) \sharp (\sharp_{15} T_{6,6})$. Then

\[
\widetilde{H}_{k}(M)=
\begin{cases}
\mathbb{Z}^{30} &k=6\\
\mathbb{Z}^{16}&n=5,7\\
\mathbb{Z}&n=12\\
0 &\text{otherwise}
\end{cases}
\]

\end{Lemma}
\begin{proof}
Let $U \subseteq T_{m,n}$ be a contractible open set. Then by the
K\"{u}nneth theorem and excision we have
\newline
\[
\widetilde{H}_{k}(T_{m,n}-U)=
\begin{cases}
\mathbb{Z} &k=m,n, \;m\ne n\\
\mathbb{Z}^{2} &k=m=n\\
0 & \text{otherwise}
\end{cases}
\]
\newline

It is clear that there exists a cofibration

\begin{equation}\label{CO}
(T_{m,n}-U)\to \sharp_{i} T_{m,n} \to    \sharp_{i} T_{m,n}/
(T_{m,n}-U)\backsimeq \sharp_{(i-1)} T_{m,n}
\end{equation}
\newline

Successive applications of the long exact sequence in homology
induced from (\ref{CO}) yields

\[
\widetilde{H}_{k}(\sharp_{16} T_{5,7})=
\begin{cases}
\mathbb{Z}^{16} &k=5,7 \\
\mathbb{Z} &k=12\\
0 & \text{otherwise}
\end{cases}
\]

and
\[
\widetilde{H}_{k}(\sharp_{15} T_{6,6})=
\begin{cases}
\mathbb{Z}^{30} &k=6 \\
\mathbb{Z} &k=12\\
0 & \text{otherwise}
\end{cases}
\]
\newline

Now we use the cofibration\\
$$((\sharp_{15}T_{6,6})-U) \to M \to M/(\sharp_{15}T_{6,6}-U) \backsimeq
\sharp_{16}T_{5,7}$$

to obtain a long exact sequence\\
$$
\cdots \to  \widetilde{H}_{k}((\sharp_{15}T_{6,6})-U) \to
\widetilde{H}_{k}(M)\to \widetilde{H}_{k}(\sharp_{16}T_{5,7}) \to
\cdots
$$
\newline

This long exact sequence in homology immediately gives
$\widetilde{H}_5(M)\cong \mathbb{Z}^{16}$ and

\begin{equation}\label{CO2}
0\to \widetilde{H}_{7}(M)\to  \mathbb{Z}^{16}\underset{j}{\to} \mathbb{Z}^{30}
\to \widetilde{H}_{6}(M)\to 0
\end{equation}

Since M is a compact, closed, orientable, manifold of dimension
$12$, then $H_{12}(M) \cong \mathbb{Z}$ and by Poincar\'{e}
duality there is an isomorphism $H^7(M) \cong H_5(M) \cong
\mathbb{Z}^{16}$. Since there is no torsion we obtain the
isomorphism $H_7(M) \cong \mathbb{Z}^{16}$ and hence the map
$j:\mathbb{Z}^{16} \rightarrow \mathbb{Z}^{30}$ must factor through
zero, giving us $H_6(M) \cong \mathbb{Z}^{30}$.

\end{proof}

Now the main result:

\begin{Theorem}
For $P = C(4,8)$, the moment angle complex $Z_P$ is not homotopy
equivalent to $(\sharp_{16} S^5 \times S^7) \sharp (\sharp_{15} S^6
\times S^6)$.
\end{Theorem}

\begin{proof}
    If $Z_P$ were homotopy equivalent to $\sharp_{16} (S^5 \times S^7) \sharp (\sharp_{15} S^6 \times
S^6)$, then by Lemma \ref{facering} and Theorem \ref{borel homotopy}
we have 

$$\pi_6(\bigvee_{16} S^5)_{(2)} \cong \pi_6((\sharp_{16} S^5
\times S^6) \sharp (\sharp_{15} S^6 \times S^6))_{(2)}$$

 We use Serre's method to calculate $\pi_{6}(M)$. Let $F$ be the homotopy fiber of
the map $M\to \prod_{16} K(\mathbb{Z},5)$.  There exists a fibration
$\prod_{16} K(\mathbb{Z},4) \to F \to M$. It is known that
$H^{*}(K(\mathbb{Z},4);\mathbb{Q})\cong \mathbb{Q}[x]$ where $|x|=4$
\cite{H}. It follows from a Serre Spectral Sequence calculation and
the universal coefficients theorem that
$\widetilde{H}_{6}(F;\mathbb{Q})\cong \mathbb{Q}^{30}$. Since $F$ is
$5$-connected, application of the Hurewicz theorem yields the
isomorphisms

$$\pi_{6}(M)\otimes \mathbb{Q} \cong \pi_{6}(F)\otimes \mathbb{Q} \cong H_{6}(F)\otimes \mathbb{Q} \cong H_{6}(F;\mathbb{Q})\cong
\mathbb{Q}^{30}$$

\vspace{.15in}

 We now can compare this result with $\pi_{6}(\vee_{16}S^{5})_{(2)}$. Using (\ref{eq}) and (\ref{f})  we know that $\pi_{\ast}(\vee_{16}S^{5}) \cong\pi_{\ast}( \prod_{16} S^{5})\times \pi _{\ast}  (\prod_{k} S^{k})$ where $k >8$.
 It suffices to check $\pi_{6}(\prod_{16}S^{5})$. But this group is all torsion. This implies
 that the group $\pi_{6}(\vee_{16}S^{5}) \otimes \mathbb{Q}$ is the trivial group.

\end{proof}
\textbf{Acknowledgements}\\

 The first author would like to thank Taras Panov
who pointed out during the International Conference on Toric
Topology in Osaka Japan, 2006 that the results in \cite{A} might be
applicable with respect to this conjecture of Bosio and Meersseman.
Both authors thank the referee for many helpful suggestions that
improved the exposition of the paper.

\end{document}